\documentclass{article} 
\usepackage{amsmath, amsfonts, amsthm, amssymb}
\usepackage{nips14submit_e,times}
\usepackage{hyperref}
\usepackage{url}

\usepackage{natbib}
\renewcommand{\bibsection}{}

\title{Exponential Concentration of a Density Functional Estimator}

\author{
Shashank Singh \\
Statistics \& Machine Learning Departments \\
Carnegie Mellon University \\
Pittsburgh, PA 15213 \\
\texttt{sss1@andrew.cmu.edu} \\
\And
Barnab\'as P\'oczos \\
Machine Learning Department \\
Carnegie Mellon University \\
Pittsburgh, PA 15213 \\
\texttt{bapoczos@cs.cmu.edu} \\
}

\renewcommand{\qed}{\quad \ensuremath{\blacksquare}}    
\newcommand{\inv}{^{-1}}                            
\newcommand{\N}{\mathbb{N}}                         
\newcommand{\R}{\mathbb{R}}                         
\newcommand{\e}{\varepsilon}                        
\newcommand{\X}{\mathcal{X}}                        
\newcommand{\Y}{\mathcal{Y}}                        
\newcommand{\Z}{\mathcal{Z}}                        
\newcommand{\E}{\mathbb{E}}                         
\newcommand{\V}{\mathbb{V}}                         
\newcommand{\pr}{\mathbb{P}}                        
\newcommand{\vi}{{\vec{i}}}                         
\newcommand{\dist}{\operatorname{dist}}             
\newcommand{\acro}[1]{\textsc{\MakeLowercase{#1}}}

\nipsfinalcopy 

\begin{document}

\maketitle

\begin{abstract}
We analyze a plug-in estimator for a large class of integral functionals of one
or more continuous probability densities. This class includes important
families of entropy, divergence, mutual information, and their conditional
versions. For densities on the $d$-dimensional unit cube $[0,1]^d$ that lie in
a $\beta$-H\"older smoothness class, we prove our estimator converges at the
rate $O \left( n^{-\frac{\beta}{\beta + d}} \right)$. Furthermore, we prove
the estimator is exponentially concentrated about its mean, whereas most
previous related results have proven only expected error bounds on estimators.
\end{abstract}

\section{Introduction}
Many important quantities in machine learning and statistics can be viewed as
integral functionals of one of more continuous probability densities; that is,
quanitities of the form
\[F(p_1,\cdots,p_k)
    = \int_{\X_1 \times \cdots \times X_k}
        f(p_1(x_1),\dots,p_k(x_k)) \, d(x_1,\dots,x_k),\]
where $p_1,\cdots,p_k$ are probability densities of random variables taking
values in $\X_1,\cdots,X_k$, respectively, and $f : \R^k \to \R$ is some
measurable function. For simplicity, we refer to such integral functionals of
densities as `density functionals'. In this paper, we study the problem of
estimating density functionals. In our framework, we assume that the underlying
distributions are not given explicitly. Only samples of $n$ independent and
identically distributed (i.i.d.) points from each of the unknown, continuous,
nonparametric distributions $p_1,\cdots,p_k$ are given.

\subsection{Motivations and Goals}
One density functional of interest is Conditional Mutual Information (CMI), a
measure of conditional dependence of random variables, which comes in several
varieties including R\'enyi-$\alpha$ and Tsallis-$\alpha$ CMI (of which Shannon
CMI is the $\alpha \to 1$ limit case). Estimating conditional dependence in a
consistent manner is a crucial problem in machine learning and statistics; for
many applications, it is important to determine how the relationship between
two variables changes when we observe additional variables. For example, upon
observing a third variable, two correlated variables may become independent,
and, similarly, two independent variables may become dependent. Hence, CMI
estimators can be used in many scientific areas to detect confounding variables
and avoid infering causation from apparent correlation
\citep{pearl98counfounding,montgomery05experiments}. Conditional dependencies
are also central to Bayesian network learning
\citep{koller09probgraphmodels,zhang11independencetest}, where CMI estimation
can be used to verify compatibility of a particular Bayes net with observed
data under a local Markov assumption.

Other important density functionals are divergences between probability
distributions, including R\'enyi-$\alpha$ \citep{renyi70probability} and
Tsallis-$\alpha$ \cite{villmann10mathematical} divergences (of which
Kullback-Leibler (\acro{KL}) divergence \citep{kullback51KL} is the
$\alpha \to 1$ limit case), and $L_p$ divergence. Divergence estimators can be
used to extend machine learning algorithms for regression, classification, and
clustering from the standard setting where inputs are finite-dimensional
feature vectors to settings where inputs are sets or distributions
\citep{poczos12kernelimages,oliva13ICML}.
Entropy and mutual information (MI) can be estimated as special cases of
divergences. Entropy estimators are used in goodness-of-fit testing
\citep{goria05new},
parameter estimation in semi-parametric models \citep{Wolsztynski85minimum},
and texture classification \citep{hero2002aes}, and MI estimators are used in
feature selection \citep{peng05feature}, clustering
\citep{aghagolzadeh07hierarchical},
optimal experimental design \citep{lewi07realtime}, and boosting and facial
expression recognition \citep{Shan05conditionalmutual}. Both entropy and mutual
information estimators are used in independent component and subspace analysis
\citep{radical03,szabo07undercomplete_TCC} and image registration
\citep{hero2002aes}. Further applications of divergence estimation are in
\citep{Leonenko-Pronzato-Savani2008}.

Despite the practical utility of density functional estimators, little is known
about their statistical performance, especially for functionals of more than
one density. In particular, few density functional estimators have known
convergence rates, and, to the best of our knowledge, no finite sample
exponential concentration bounds have been derived for general density
functional estimators. One consequence of this exponential bound is that, using
a union bound, we can guarantee accuracy of multiple estimates
simultaneously. For example, \citep{liu12exponential} shows how this can be
applied to optimally analyze forest density estimation algorithms. Because the
CMI of variables $X$ and $Y$ given a third variable $Z$ is zero if and only $X$
and $Y$ are conditionally independent given $Z$, by estimating CMI with a
confidence interval, we can test for conditional independence with bounded type
I error probabilty.

{\bf Our main contribution} is to derive convergence rates and an exponential
concentration inequality for a particular, consistent, nonparametric estimator
for large class of density functionals, including conditional density
functionals. We also apply our concentration inequality to the important case
of R\'enyi-$\alpha$ CMI.

\subsection{Related Work}
Although lower bounds are not known for estimation of general density
functionals (of arbitrarily many densities), \citep{birge95estimation} lower
bounded the convergence rate for estimators of functionals of a single density
(e.g., entropy functionals) by $O\left( n^{-4\beta/(4\beta + d)} \right)$.
\citep{krishnamurthy14divergences} extended
this lower bound to the two-density cases of $L_2$, R\'enyi-$\alpha$, and
Tsallis-$\alpha$ divergences and gave plug-in estimators which achieve this
rate. These estimators enjoy the parametric rate of $O\left( n^{-1/2} \right)$
when $\beta > d/4$, and work by optimally estimating the density and then
applying a correction to the plug-in estimate. In contrast, our estimator
undersmooths the density, and converges at a slower rate of
$O\left( n^{-\beta/(\beta + d)} \right)$ when $\beta < d$ (and the parametric
rate $O\left( n^{-1/2} \right)$ when $\beta \geq d$), but obeys an exponential
concentration inequality, which is not known for the estimators of
\citep{krishnamurthy14divergences}.

Another exception for $f$-divergences is provided by
\citep{nguyen10estimating}, using empirical risk minimization. This approach involves
solving an $\infty$-dimensional convex minimization problem which be reduced to
an $n$-dimensional problem for certain function classes defined by reproducing
kernel Hilbert spaces ($n$ is the sample size). When $n$ is large, these
optimization problems can still be very demanding. They studied the estimator's
convergence rate, but did not derive concentration bounds.

A number of papers have studied $k$-nearest-neighbors estimators, primarily for
R\'enyi$\alpha$ density functionals including entropy \citep{leonenko08tatra},
divergence \citep{Wang-Kulkarni-Verdu2009} and conditional divergence and MI
\citep{poczos12AISTATS}. These estimators work directly, without the
intermediate density estimation step, and generally have proofs of consistency,
but their convergence rates and dependence on $k$, $\alpha$, and the dimension
are unknown. One recent exception is a $k$-nearest-neighbors
based estimator that converges at the parametric rate when $\beta > d$, using
an optimally weighted ensemble of weak estimators
\citep{sricharan13ensemble, moon14ensemble}. These estimators appear to perform
well in higher dimensions, but rates for these estimators require that
$k \to \infty$ as $n \to \infty$, causing computational difficulties for large
samples.

Although the literature on dependence measures is huge, few estimators have
been generalized to the conditional case
\citep{fukumizu08kernelmeasures,reddi13dependence}. There is some work on
testing conditional dependence
\citep{su08dependencetest,bouezmarni09dependencetest}, but, unlike CMI
estimation, these tests are intended to simply accept or reject the hypothesis
that variables are conditionally independent, rather than to measure
conditional dependence. Our exponential concentration inequality also suggests
a new test for conditional independence.

This paper continues a line of work begun in \citep{liu12exponential} and
continued in \citep{singh14exponential}.
\citep{liu12exponential} proved an exponential concentration inequality for an
estimator of Shannon entropy and MI in the $2$-dimensional case.
\citep{singh14exponential} used similar techniques to derive an exponential
concentration inequality for an estimator of R\'enyi-$\alpha$ divergence in
$d$ dimensions, for a larger family of densities. Both used plug-in estimators
based on a mirrored kernel density estimator (KDE) on $[0,1]^d$. Our work
generalizes these results to a much larger class of density functionals, as
well as to conditional density functionals
(see Section 6). In particular, we use a plug-in estimator for general density
functionals based on the same mirrored KDE, and also use some lemmas regarding
this KDE proven in \citep{singh14exponential}. By considering the more general
density functional case, we are also able to significantly simplify the proofs
of the convergence rate and exponential concentration inequality.

\subsection*{Organization}
In Section 2, we establish the theoretical context of our work, including
notation, the precise problem statement, and our estimator. In Section 3, we
outline our main theoretical results and state some consequences. Sections 4
and 5 give precise statements and proofs of the results in Section 3. Finally,
in Section 6, we extend our results to conditional density functionals,
and state the consequences in the particular case of R\'enyi-$\alpha$ CMI.

\section{Density Functional Estimator}
\subsection{Notation}
For an integer $k$, $[k] = \{1,\cdots,k\}$ denotes the set of positive integers
at most $k$. Using the notation of multi-indices common in multivariate
calculus, $\N^d$ denotes the set of $d$-tuples of non-negative integers, which
we denote with a vector symbol $\vec\cdot$, and, for $\vi \in \N^d$,
\[D^\vi := \frac{\partial^{|\vi|}}{\partial^{i_1}x_1\cdots\partial^{i_d}x_d}
    \quad \mbox{ and } \quad
    |\vi| = \sum_{k = 1}^d i_k.
\]
For fixed $\beta, L > 0$, $r \geq 1$, and a positive integer $d$, we will work
with densities in the following bounded subset of a $\beta$-H\"older space:
\begin{align}
C_{L,r}^\beta([0,1]^d)
    := \left\{ p : [0,1]^d \to \R \middle|
            \sup_{\substack{x \neq y \in D\\|\vi| = \ell}}
            \frac{|D^\vi p(x) - D^\vi p(y)|}{\|x - y\|^{(\beta - \ell)}}
    \right\},
\label{eq:holder}
\end{align}
where $\ell = \lfloor \beta \rfloor$ is the greatest integer \emph{strictly}
less than $\beta$, and $\|\cdot\|_r : \R^d \to \R$ is the usual $r$-norm. To
correct for boundary bias, we will require the densities to be nearly constant
near the boundary of $[0,1]^d$, in that their derivatives vanish at the
boundary. Hence, we work with densities in
\begin{align}
\Sigma(\beta,L,r,d)
    := \left\{ p \in C_{L,r}^\beta([0,1]^d) \middle|
            \max_{1 \leq |\vi| \leq \ell} |D^\vi p(x)| \to 0
            \mbox{ as }
            \dist(x,\partial[0,1]^d) \to 0
    \right\},
\label{eq:bdd_holder}
\end{align}
where
$\partial[0,1]^d = \{x \in [0,1]^d : x_j \in \{0,1\}
    \mbox{ for some } j \in [d]\}$.

\subsection{Problem Statement}
For each $i \in [k]$ let $X_i$ be a $d_i$-dimensional random vector taking
values in $\X_i := [0,1]^{d_i}$, distributed according to a density
$p_i : \X \to \R$. For an appropriately smooth function $f : \R^k \to \R$, we
are interested in a using random sample of $n$ i.i.d. points from the
distribution of each $X_i$ to estimate
\begin{equation}
F(p_1,\dots,p_k)
    := \int_{\X_1 \times \cdots \times \X_k}
        f(p_1(x_1),\dots,p_k(x_k)) \, d(x_1,\dots,x_k).
\label{eq:df_form}
\end{equation}

\subsection{Estimator}
For a fixed bandwidth $h$, we first use the mirrored kernel density estimator
(KDE) $\hat p_i$ described in \citep{singh14exponential} to estimate each
density $p_i$. We then use a plug-in estimate of $F(p_1,\dots,p_k)$.
\[F(\hat p_1,\dots,\hat p_k)
    := \int_{\X_1 \times \cdots \times \X_k}
        f(\hat p_1(x_1),\dots,\hat p_k(x_k)) \, d(x_1,\dots,x_k).
\]
Our main results generalize those of \citep{singh14exponential} to a broader
class of density functionals.

\section{Main Results}
In this section, we outline our main theoretical results, proven in Sections 4
and 5, and also discuss some important corollaries.

We decompose the estimatator's error into bias and a variance-like terms via
the triangle inequality:
\begin{align*}
|F(\hat p_1,\dots,\hat p_k) - F(p_1,\dots,p_k)|
 &  \leq \underbrace{|F(\hat p_1,\dots,\hat p_k)
                                            - \E F(\hat p_1,\dots,\hat p_k)|}
        _{\text{variance-like term}}    \\
 &  + \underbrace{|\E F(\hat p_1,\dots,\hat p_k) - F(p_1,\dots,p_k)|}
        _{\text{bias term}}.
\end{align*}
We will prove the ``variance'' bound
\begin{equation}
\pr \left( |F(\hat p_1,\dots,\hat p_k)
                                - \E F(\hat p_1,\dots,\hat p_k)| > \e \right)
    \leq 2\exp \left( -\frac{2\e^2n}{C_V^2} \right)
\label{ineq:var_bdd}
\end{equation}
for all $\e > 0$ and the bias bound
\begin{equation}
|\E F(\hat p_1,\dots,\hat p_k) - F(p_1,\dots,p_k)|
    \leq C_B \left( h^\beta + h^{2\beta} + \frac{1}{nh^d} \right),
\label{ineq:bias_bdd}
\end{equation}
where $d := \max_i d_i$, and $C_V$ and $C_B$ are constant in the sample size
$n$ and bandwidth $h$ for exact values. To the best of our knowledge, this is
the first time an exponential inequality like (\ref{ineq:var_bdd}) has been
established for general density functional estimation. This variance bound does
not depend on $h$ and the bias bound is minimized by
$h \asymp n^{-\frac{1}{\beta + d}}$, we have the convergence rate
\[|\E F(\hat p_1,\dots,\hat p_k) - F(p_1,\dots,p_k)|
    \in O \left( n^{-\frac{\beta}{\beta + d}} \right).
\]
It is interesting to note that, in optimizing the bandwidth for our density
functional estimate, we use a smaller bandwidth than is optimal for minimizing
the bias of the KDE. Intuitively, this reflects the fact that the plug-in
estimator, as an integral functional, performs some additional smoothing.

We can use our exponential concentration bound to obtain a bound on the true
variance of $F(\hat p_1,\dots,\hat p_k)$. If $G : [0,\infty) \to \R$ denotes
the cumulative distribution function of the squared deviation of
$F(\hat p_1,\dots,\hat p_k)$ from its mean, then
\[1 - G(\e)
    = \pr \left( \left( F(\hat p_1,\dots,\hat p_k)
                        - \E F(\hat p_1,\dots,\hat p_k) \right)^2 > \e \right)
    \leq 2\exp \left( -\frac{2\e n}{C_V^2} \right).
\]
Thus,
\begin{align*}
\V[F(\hat p_1,\dots,\hat p_k)]
 &  = \E\left[ \left( F(\hat p_1,\dots,\hat p_k)
                        - \E F(\hat p_1,\dots,\hat p_k) \right)^2 \right]   \\
 &  = \int_0^\infty 1 - G(\e) \, d\e
    \leq 2\int_0^\infty \exp \left( -\frac{2\e n}{C_V^2} \right)
    = C_V^2n\inv.
\end{align*}
We then have a mean squared error of
\[\E \left[ \left( F(\hat p_1,\dots,\hat p_k)
                                        - F(p_1,\dots,p_k) \right)^2 \right]
    \in O \left( n\inv + n^{-\frac{2\beta}{\beta + d}} \right),
\]
which is in $O(n\inv)$ if $\beta \geq d$ and
$O \left( n^{-\frac{2\beta}{\beta + d}} \right)$ otherwise.

It should be noted that the constants in both the bias bound and the variance
bound depend exponentially on the dimension $d$. Lower bounds in terms of $d$
are unknown for estimating most density functionals of interest, and an
important open problem is whether this dependence can be made asymptotically
better than exponential.

\section{Bias Bound}
In this section, we precisely state and prove the bound on the bias of our
density functional estimator, as introduced in Section 3.

Assume each $p_i \in \Sigma(\beta,L,r,d)$ (for $i \in [k]$), assume
$f : \R^k \to \R$ is twice continuously differentiable, with first and second
derivatives all bounded in magnitude by some $C_f \in \R$,
\footnote{If $p_1(\X_1)\times\cdots\times p_k(\X_k)$ is known to lie within
some cube $[\kappa_1,\kappa_2]^k$, then it suffices for $f$ to be twice
continuously differentiable on $[\kappa_1,\kappa_2]^k$ (and the boundedness
condition follows immediately). This will be important for our application to
R\'enyi-$\alpha$ Conditional Mutual Information.}
and assume the kernel $K : \R \to \R$ has bounded support $[-1,1]$ and
satisfies
\[\int_{-1}^1 K(u) \, du = 1
    \quad\mbox{ and }\quad
    \int_{-1}^1 u^j K(u) \, du = 0
    \quad\mbox{ for all } j \in \{1,\dots,\ell\}.
\]
Then, there exists a constant $C_B \in \R$ such that
\[|\E F(\hat p_1,\dots,\hat p_k) - F(p_1,\dots,p_k)|
    \leq C_B \left( h^\beta + h^{2\beta} + \frac{1}{nh^d} \right).
\]

\subsection{Proof of Bias Bound}
By Taylor's Theorem,
$\forall x = (x_1,\dots,x_k) \in \X_1\times\cdots\times\X_k$, for some
$\xi \in \R^k$ on the line segment between
$\hat p(x) := (\hat p_1(x_1),\dots,\hat p_k(x_k))$ and
$p(x) := (p_1(x_1),\dots,p_k(x_k))$, letting $H_f$ denote the Hessian of $f$
\begin{align*}
|\E f(\hat p(x)) - f(p(x))|
 &  = \left| \E (\nabla f)(p(x))
        \cdot \left( \hat p(x) - p(x) \right)
    + \frac12 (\hat p(x) - p(x))^TH_f(\xi)(\hat p(x) - p(x)) \right|    \\
 &  \hspace{-5mm}
    \leq C_f \left( \sum_{i = 1}^k |B_{p_i}(x_i)|
                    + \sum_{i < j \leq k} |B_{p_i}(x_i)B_{p_j}(x_j)|
                    + \sum_{i = 1}^k \E[\hat p_i(x_i) - p_i(x_i)]^2
    \right)
\end{align*}
where we used that $\hat p_i$ and $\hat p_j$ are independent for $i \neq j$.
Applying H\"older's Inequality,
\begin{align*}
 & |\E F(\hat p_1,\dots,\hat p_k) - F(p_1,\dots,p_k)|
    \leq \int_{\X_1\times\cdots\times\X_k} |\E f(\hat p(x)) - f(p(x))| \, dx \\
 &  \leq C_f\hspace{-1mm}
            \left( \sum_{i = 1}^k \int_{\X_i} \hspace{-1.5mm}|B_{p_i}(x_i)|
                                    + \E[\hat p_i(x_i) - p_i(x_i)]^2 \, dx_i
                    + \hspace{-2mm}\sum_{i < j \leq k}
                    \int_{\X_i} \hspace{-1.5mm}|B_{p_i}(x_i)| \, dx_i
                    \int_{\X_j} \hspace{-1.5mm}|B_{p_j}(x_j)| \, dx_j
            \right) \\
 &  \leq C_f
            \Bigg( \sum_{i = 1}^k \sqrt{\int_{\X_i} B_{p_i}^2(x_i) \, dx_i}
                        + \int_{\X_i} \E[\hat p_i(x_i) - p_i(x_i)]^2 \, dx_i \\
 & \hspace{70mm}    + \sum_{i < j \leq k}
                    \sqrt{\int_{\X_i} B_{p_i}^2(x_i) \, dx_i
                    \int_{\X_j} B_{p_j}^2(x_j) \, dx_j}
            \Bigg).
\end{align*}

We now make use of the so-called Bias Lemma proven by
\citep{singh14exponential}, which bounds the integrated squared bias of the
mirrored KDE $\hat p$ on $[0,1]^d$ for an arbitrary
$p \in \Sigma(\beta,L,r,d)$. Writing the bias of $\hat p$ at $x \in [0,1]^d$ as
$B_p(x) = \E\hat p(x) - p(x)$, \citep{singh14exponential} showed that there
exists $C > 0$ constant in $n$ and $h$ such that
\begin{equation}
\int_{[0,1]^d} B_p^2(x) \, dx \leq Ch^{2\beta}.
\label{ineq:bias_lemma}
\end{equation}
Applying the Bias Lemma and certain standard results in kernel density
estimation (see, for example, Propositions 1.1 and 1.2 of
\citep{Tsybakov:2008:INE:1522486}) gives
\begin{align}
\notag
|\E F(\hat p_1,\dots,\hat p_k) - F(p_1,\dots,p_k)|
 &  \leq C \left( k^2h^\beta + kh^{2\beta} \right) + \frac{\|K\|_1^d}{nh^d}
 &  \leq C_B \left( h^\beta + h^{2\beta} + \frac{1}{nh^d} \right),
\end{align}
where $\|K\|_1$ denotes the $1$-norm of the kernel. \qed

\section{Variance Bound}
In this section, we precisely state and prove the exponential concentration
inequality for our density functional estimator, as introduced in Section 3.
Assume that $f$ is Lipschitz continuous with constant $C_f$ in the $1$-norm
on $p_1(\X_1) \times\cdots\times p_k(\X_k)$ (i.e.,
\begin{equation}
|f(x) - f(y)| \leq C_f\sum_{k = 1}^\infty |x_i - y_i|,
    \quad \forall x,y \in p_1(\X_1) \times\cdots\times p_k(\X_k)).
\label{ineq:lip_cond}
\end{equation}

and assume the kernel $K \in L_1(\R)$ (i.e., it has finite $1$-norm). Then,
there exists a constant $C_V \in \R$ such that $\forall \e > 0$,
\[\pr\left(|F(\hat p_1,\dots,\hat p_k) - \E F(\hat p_1,\dots,\hat p_k)|\right)
    \leq 2\exp \left( -\frac{2\e^2n}{C_V^2} \right).\]
Note that, while we require no assumptions on the densities here, in certain
specific applications, such us for some R\'enyi-$\alpha$ quantities, where
$f = \log$, assumptions such as lower bounds on the density may be needed to
ensure $f$ is Lipschitz on its domain.

\subsection{Proof of Variance Bound}
Consider i.i.d. samples $(x^1_1,\dots,x^n_k) \in \X_1\times\cdots\times\X_k$
drawn according to the product distribution $p = p_1 \times \cdots p_k$. In
anticipation of using McDiarmid's Inequality \citep{McDiarmid1989}, let
$\hat p_j'$ denote the $j^{th}$ mirrored KDE when the sample $x^i_j$ is
replaced by new sample $(x^i_j)'$. Then, applying the Lipschitz condition
(\ref{ineq:lip_cond}) on $f$,
\[|F(\hat p_1,\dots,\hat p_k) - F(\hat p_1,\dots,\hat p_j',\dots,\hat p_k)|
    \leq C_f\int_{\X_j} |p_j(x) - p_j'(x)| \, dx,
\]
since most terms of the sum in (\ref{ineq:lip_cond}) are zero. Expanding the
definition of the kernel density estimates $\hat p_j$ and $\hat p_j'$ and
noting that most terms of the mirrored KDEs $\hat p_j$ and $\hat p_j'$ are
identical gives
\begin{align*}
|F(\hat p_1,\dots,\hat p_k) - F(\hat p_1,\dots,\hat p_j',\dots,\hat p_k)|
 &  = \frac{C_f}{nh^{d_j}} \int_{\X_j} \left|
        K_{d_j}\left( \frac{x - x^i_j}{h} \right)
            - K_{d_j}\left( \frac{x - (x^i_j)'}{h} \right) \right| \, dx
\end{align*}
where $K_{d_j}$ denotes the $d_j$-dimensional mirrored product kernel based on
$K$. Performing a change of variables to remove $h$ and applying the triangle
inequality followed by the bound on the integral of the mirrored kernel proven
in \citep{singh14exponential},
\begin{align}
\notag
|F(\hat p_1,\dots,\hat p_k) - F(\hat p_1,\dots,\hat p_j',\dots,\hat p_k)|
 &  \leq \frac{C_f}{n} \int_{h\X_j}
            \left| K_{d_j}(x - x^i_j) - K_{d_j}(x - (x^i_j)')\right| \, dx  \\
 &  \leq \frac{2C_f}{n} \int_{[-1,1]^{d_j}} \hspace{-5mm}|K_{d_j}(x)| \, dx
    \leq \frac{2C_f}{n} \|K\|_1^{d_j}
    = \frac{C_V}{n},
\label{ineq:mcd_cond}
\end{align}
for $C_V = 2C_f \max_j \|K\|_1^{d_j}$. Since $F(\hat p_1,\dots,\hat p_k)$
depends on $kn$ independent variables, McDiarmid's Inequality then gives, for
any $\e > 0$,
\[\mathbb{P}\left( |F(\hat p_1,\dots,\hat p_k) - F(p_1,\dots,p_k)| > \e \right)
    \leq 2\exp\left( -\frac{2\e^2}{knC_V^2/n^2} \right)
    = 2\exp\left( -\frac{2\e^2n}{kC_V^2} \right). \qed
\]

\section{Extension to Conditional Density Functionals}
Our convergence result and concentration bound can be fairly easily adapted to
to KDE-based plug-in estimators for many functionals of interest, including
R\'enyi-$\alpha$ and Tsallis-$\alpha$ entropy, divergence, and MI, and
$L_p$ norms and distances, which have either the same or analytically similar
forms as as the functional (\ref{eq:df_form}). As long as the density of the
variable being conditioned on is lower bounded on its domain, our results also
extend to conditional density functionals of the form
\footnote{We abuse notation slightly and also use $P$ to denote all of its
marginal densities.}
\begin{equation}
F(P)
    = \int_\Z P(z) f\left(
            \int_{\X_1 \times \cdots \times \X_k} g\left(
                \frac{P(x_1,z)}{P(z)},
                \frac{P(x_2,z)}{P(z)},
                \dots,
                \frac{P(x_k,z)}{P(z)}
            \right) \, d(x_1,\dots,x_k)
        \right) \, dz
\label{eq:cond_df_form}
\end{equation}
including, for example, R\'enyi-$\alpha$ conditional entropy, divergence, and
mutual information, where $f$ is the function
$x \mapsto \frac{1}{1 - \alpha} \log(x)$. The proof of this extension for
general $k$ is essentially the same as for the case $k = 1$, and so, for
notational simplicity, we demonstrate the latter.

\subsection{Problem Statement, Assumptions, and Estimator}
For given dimensions $d_x,d_z \geq 1$, consider random vectors $X$ and $Z$
distributed on unit cubes $\X := [0,1]^{d_x}$ and $\Z := [0,1]^{d_z}$
according to a joint density $P : \X \times \Z \to \R$. We use a random sample
of $2n$ i.i.d. points from $P$ to estimate a conditional density functional
$F(P)$, where $F$ has the form (\ref{eq:cond_df_form}). 

Suppose that $P$ is in the H\"older class $\Sigma(\beta,L,r,d_x + d_z)$, noting
that this implies an analogous condition on each marginal of $P$, and suppose
that $P$ bounded below and above, i.e.,
$0 < \kappa_1 := \inf_{x \in \X,z \in \Z} P(z)$ and
$\infty > \kappa_2 := \inf_{x \in \X,z \in \Z} P(x,z)$. Suppose also that $f$
and $g$ are continuously differentiable, with
\begin{equation}
C_f := \sup_{x \in [c_g,C_g]} |f(x)|
    \quad \mbox{ and } \quad
    C_{f'} := \sup_{x \in [c_g,C_g]} |f'(x)|,
\label{ineq:f_bounds}
\end{equation}
where
\[c_g := \inf g\left(\left[0, \frac{\kappa_2}{\kappa_1}\right]\right)
    \quad \mbox{ and } \quad
    C_g := \sup g\left(\left[0, \frac{\kappa_2}{\kappa_1}\right]\right).
\]
After estimating the densities $P(z)$ and $P(x,z)$ by their mirrored KDEs,
using $n$ independent data samples for each, we clip the estimates of $P(x,z)$
and $P(z)$ below by $\kappa_1$ and above by $\kappa_2$ and denote the resulting
density estimates by $\hat P$. Our estimate $F(\hat P)$ for $F(P)$ is simply
the result of plugging $\hat P$ into equation (\ref{eq:cond_df_form}).

\subsection{Proof of Bounds for Conditional Density Functionals}
We bound the error of $F(\hat P)$ in terms of the error of estimating the
corresponding unconditional density functional using our previous estimator,
and then apply our previous results.

Suppose $P_1$ is either the true density $P$ or a plug-in estimate of $P$
computed as described above, and $P_2$ is a plug-in estimate of $P$ computed in
the same manner but using a different data sample.
Applying the triangle inequality twice,
\begin{align*}
|F(P_1) - F(P_2)|
 &  \leq
    \int_\Z
        \left| P_1(z) f\left(
                    \int_\X g\left(\frac{P_1(x,z)}{P_1(z)} \right) \, dx \right)
            - P_2(z) f\left(
                    \int_\X g\left(\frac{P_1(x,z)}{P_1(z)} \right) \, dx \right)
        \right|
   \\
 &  +
        \left| P_2(z) f\left(
                    \int_\X g\left(\frac{P_1(x,z)}{P_1(z)} \right) \, dx \right)
            - P_2(z) f\left(
                    \int_\X g\left(\frac{P_2(x,z)}{P_2(z)} \right) \, dx \right)
        \right| \, dz   \\
 &  \leq
    \int_\Z
        |P_1(z) - P_2(z)|
        \left|
            f\left( \int_\X g\left(\frac{P_1(x,z)}{P_1(z)} \right) \, dx \right)
        \right|
   \\
 &  + P_2(z)
        \left| f\left(
                    \int_\X g\left(\frac{P_1(x,z)}{P_1(z)} \right) \, dx \right)
            - f\left(
                    \int_\X g\left(\frac{P_2(x,z)}{P_2(z)} \right) \, dx \right)
        \right| \, dz
\end{align*}
Applying the Mean Value Theorem and the bounds in (\ref{ineq:f_bounds}) gives
\begin{align*}
|F(P_1) - F(P_2)|
 &  \leq \hspace{-1mm}\int_\Z\hspace{-1mm}C_f|P_1(z) - P_2(z)|
    + \kappa_2C_{f'}
        \left| \int_\X \hspace{-1mm}g\left(\frac{P_1(x,z)}{P_1(z)} \right)
            - g\left(\frac{P_2(x,z)}{P_2(z)} \right) dx
        \right| dz  \\
 &  = \hspace{-1mm}\int_\Z\hspace{-1mm}C_f |P_1(z) - P_2(z)|
    + \kappa_2C_{f'} \left| G_{P_1(z)}(P_1(\cdot,z)) - G_{P_2(z)}(P_2(\cdot,z)) \right| dz,
\end{align*}
where $G_z$ is the density functional
\[G_{P(z)}(Q) = \int_\X g\left( \frac{Q(x)}{P(z)} \right) \, dx.\]
Note that, since the data are split to estimate $P(z)$ and $P(x,z)$,
$G_{\hat P(z)}(\hat P(\cdot,z))$ depends on each data point through only one of
these KDEs. In the case that $P_1$ is the true density $P$, taking the
expectation and using Fubini's Theorem gives
\begin{align*}
\E|F(P) - F(\hat P)|
 &  \leq \hspace{-1mm}\int_\Z\hspace{-1mm}C_f \E |P(z) - \hat P(z)|
    + \kappa_2 C_{f'} \E \left| G_{P(z)}(P(\cdot,z)) - G_{\hat P(z)}(\hat
P(\cdot,z)) \right| dz, \\
 &  \leq \hspace{-1mm}C_f\sqrt{\int_\Z\hspace{-1mm} \E (P(z) - \hat P(z))^2 dz}
    + 2\kappa_2 C_{f'} C_B \left( h^\beta + h^{2\beta} + \frac{1}{nh^d} \right) \\
 &  \leq (2\kappa_2 C_{f'} C_B + C_fC)\left( h^\beta + h^{2\beta} + \frac{1}{nh^d} \right)
\end{align*}
applying H\"older's Inequality and our bias bound (\ref{ineq:bias_bdd}),
followed by the bias lemma (\ref{ineq:bias_lemma}). This extends our bias bound
to conditional density functionals. For the variance bound, consider the case
where $P_1$ and $P_2$ are each mirrored KDE estimates of $P$, but with one data
point resampled (as in the proof of the variance bound, setting up to use
McDiarmid's Inequality). By the same sequence of steps used to show
(\ref{ineq:mcd_cond}),
\[\int_\Z |P_1(z) - P_2(z)| \, dz \leq \frac{2\|K\|_1^{d_z}}{n},\]
and
\begin{align*}
\int_\Z \left| G_{P(z)}(P(\cdot,z)) - G_{\hat P(z)}(\hat P(\cdot,z)) \right| dz
    \leq \frac{C_V}{n}.
\end{align*}
(by casing on whether the resampled data point was used to estimate $P(x,z)$ or
$P(z)$), for an appropriate $C_V$ depending on
$\sup_{x \in [\kappa_1/\kappa_2,\kappa_2/\kappa_1]} |g'(x)|$.
Then, by McDiarmid's Inequality,
\[\mathbb{P}\left( |F(\hat p_1,\dots,\hat p_k) - F(p_1,\dots,p_k)| > \e \right)
    = 2\exp\left( -\frac{\e^2n}{4C_V^2} \right). \qed
\]

\subsection{Application to R\'enyi-$\alpha$ Conditional Mutual Information}
As an example, we demonstrate our concentration inequality to the R\'enyi-$\alpha$
Conditional Mutual Information (CMI). Consider random vectors $X,Y$, and $Z$ on
$\X = [0,1]^{d_x}$, $\Y = [0,1]^{d_y}$, $\Z = [0,1]^{d_z}$, respectively.
$\alpha \in (0,1) \cup (1,\infty)$, the R\'enyi-$\alpha$ CMI of $X$ and $Y$
given $Z$ is
\begin{equation}
I(X; Y | Z)
    = \frac{1}{1 - \alpha} \int_\Z P(z) \log
            \int_{\X \times \Y} \left( \frac{P(x,y,z)}{P(z)} \right)^\alpha
            \left( \frac{P(x,z)P(y,z)}{P(z)^2} \right)^{1 - \alpha} \, d(x,y)
      \, dz.
\label{eq:R_CMI}
\end{equation}
In this case, the estimator which plugs mirrored KDEs for $P(x,y,z)$, $P(x,z)$,
$P(y,z)$, and $P(z)$ into (\ref{eq:R_CMI}) obeys the concentration inequality
(\ref{ineq:var_bdd}) with $C_V = \kappa^*\|K\|_1^{d_x + d_y + d_z}$, where
$\kappa^*$ depends only on $\alpha$, $\kappa_1$, and $\kappa_2$.

\subsubsection*{References}
\setlength{\bibsep}{0.0pt}
{
\bibliographystyle{plain}
\bibliography{biblio}

\begin{thebibliography}{10}

\bibitem{aghagolzadeh07hierarchical}
M.~Aghagolzadeh, H.~Soltanian-Zadeh, B.~Araabi, and A.~Aghagolzadeh.
\newblock A hierarchical clustering based on mutual information maximization.
\newblock In {\em in Proc. of IEEE International Conference on Image
  Processing}, pages 277--280, 2007.

\bibitem{birge95estimation}
L.~Birge and P.~Massart.
\newblock Estimation of integral functions of a density.
\newblock {\em A. Statistics}, 23:11--29, 1995.

\bibitem{bouezmarni09dependencetest}
T.~Bouezmarni, J.~Rombouts, and A.~Taamouti.
\newblock A nonparametric copula based test for conditional independence with
  applications to granger causality, 2009.
\newblock Technical report, Universidad Carlos III, Departamento de Economia.

\bibitem{fukumizu08kernelmeasures}
K.~Fukumizu, A.~Gretton, X.~Sun, and B.~Schoelkopf.
\newblock Kernel measures of conditional dependence.
\newblock In {\em Neural Information Processing Systems (NIPS)}, 2008.

\bibitem{goria05new}
M.~N. Goria, N.~N. Leonenko, V.~V. Mergel, and P.~L.~Novi Inverardi.
\newblock A new class of random vector entropy estimators and its applications
  in testing statistical hypotheses.
\newblock {\em J. Nonparametric Statistics}, 17:277--297, 2005.

\bibitem{hero2002aes}
A.~O. Hero, B.~Ma, O.~J.~J. Michel, and J.~Gorman.
\newblock Applications of entropic spanning graphs.
\newblock {\em IEEE Signal Processing Magazine}, 19(5):85--95, 2002.

\bibitem{koller09probgraphmodels}
D.~Koller and N.~Friedman.
\newblock {\em Probabilistic Graphical Models: Principles and Techniques}.
\newblock MIT Press, Cambridge, MA, 2009.

\bibitem{krishnamurthy14divergences}
A.~Krishnamurthy, K.~Kandasamy, B.~Poczos, and L.~Wasserman.
\newblock Nonparametric estimation of renyi divergence and friends.
\newblock In {\em International Conference on Machine Learning (ICML)}, 2014.

\bibitem{kullback51KL}
S.~Kullback and R.A. Leibler.
\newblock On information and sufficiency.
\newblock {\em Annals of Mathematical Statistics}, 22:79--86, 1951.

\bibitem{radical03}
E.~G. Learned-Miller and J.~W. Fisher.
\newblock {ICA} using spacings estimates of entropy.
\newblock {\em J. Machine Learning Research}, 4:1271--1295, 2003.

\bibitem{Leonenko-Pronzato-Savani2008}
N.~Leonenko, L.~Pronzato, and V.~Savani.
\newblock A class of {R}{\'e}nyi information estimators for multidimensional
  densities.
\newblock {\em Annals of Statistics}, 36(5):2153--2182, 2008.

\bibitem{leonenko08tatra}
N.~Leonenko, L.~Pronzato, and V.~Savani.
\newblock Estimation of entropies and divergences via nearest neighbours.
\newblock {\em Tatra Mt. Mathematical Publications}, 39, 2008.

\bibitem{lewi07realtime}
J.~Lewi, R.~Butera, and L.~Paninski.
\newblock Real-time adaptive information-theoretic optimization of
  neurophysiology experiments.
\newblock In {\em Advances in Neural Information Processing Systems},
  volume~19, 2007.

\bibitem{liu12exponential}
H.~Liu, J.~Lafferty, and L.~Wasserman.
\newblock Exponential concentration inequality for mutual information
  estimation.
\newblock In {\em Neural Information Processing Systems (NIPS)}, 2012.

\bibitem{McDiarmid1989}
C.~McDiarmid.
\newblock On the method of bounded differences.
\newblock {\em Surveys in Combinatorics}, 141:148--188, 1989.

\bibitem{montgomery05experiments}
D.~Montgomery.
\newblock {\em Design and Analysis of Experiments}.
\newblock John Wiley and Sons, 2005.

\bibitem{moon14ensemble}
K.R. Moon and A.O. Hero.
\newblock Ensemble estimation of multivariate f-divergence.
\newblock In {\em Information Theory (ISIT), 2014 IEEE International Symposium
  on}, pages 356--360, June 2014.

\bibitem{nguyen10estimating}
X.~Nguyen, M.J. Wainwright, and M.I. Jordan.
\newblock Estimating divergence functionals and the likelihood ratio by convex
  risk minimization.
\newblock {\em IEEE Trans. on Information Theory.}, 2010.

\bibitem{oliva13ICML}
J.~Oliva, B.~Poczos, and J.~Schneider.
\newblock Distribution to distribution regression.
\newblock In {\em International Conference on Machine Learning (ICML)}, 2013.

\bibitem{pearl98counfounding}
J.~Pearl.
\newblock Why there is no statistical test for confounding, why many think
  there is, and why they are almost right, 1998.
\newblock UCLA Computer Science Department Technical Report R-256.

\bibitem{peng05feature}
H.~Peng and C.~Dind.
\newblock Feature selection based on mutual information: Criteria of
  max-dependency, max-relevance, and min-redundancy.
\newblock {\em IEEE Trans On Pattern Analysis and Machine Intelligence}, 27,
  2005.

\bibitem{poczos12AISTATS}
B.~Poczos and J.~Schneider.
\newblock Nonparametric estimation of conditional information and divergences.
\newblock In {\em International Conference on AI and Statistics (AISTATS)},
  volume~20 of {\em JMLR Workshop and Conference Proceedings}, 2012.

\bibitem{poczos12kernelimages}
B.~Poczos, L.~Xiong, D.~Sutherland, and J.~Schneider.
\newblock Nonparametric kernel estimators for image classification.
\newblock In {\em 25th IEEE Conference on Computer Vision and Pattern
  Recognition (CVPR)}, 2012.

\bibitem{reddi13dependence}
S.~J. Reddi and B.~Poczos.
\newblock Scale invariant conditional dependence measures.
\newblock In {\em International Conference on Machine Learning (ICML)}, 2013.

\bibitem{renyi70probability}
A.~R\'enyi.
\newblock {\em Probability Theory}.
\newblock North-Holland Publishing Company, Amsterdam, 1970.

\bibitem{Shan05conditionalmutual}
C.~Shan, S.~Gong, and P.~W. Mcowan.
\newblock Conditional mutual information based boosting for facial expression
  recognition.
\newblock In {\em British Machine Vision Conference (BMVC)}, 2005.

\bibitem{singh14exponential}
S.~Singh and B.~Poczos.
\newblock Generalized exponential concentration inequality for r\'enyi
  divergence estimation.
\newblock In {\em International Conference on Machine Learning (ICML)}, 2014.

\bibitem{sricharan13ensemble}
K.~Sricharan, D.~Wei, and A.~Hero.
\newblock Ensemble estimators for multivariate entropy estimation, 2013.

\bibitem{su08dependencetest}
L.~Su and H.~White.
\newblock A nonparametric {H}ellinger metric test for conditional independence.
\newblock {\em Econometric Theory}, 24:829--864, 2008.

\bibitem{szabo07undercomplete_TCC}
Z.~Szab{\'o}, B.~P{\'o}czos, and A.~L{\H{o}}rincz.
\newblock Undercomplete blind subspace deconvolution.
\newblock {\em J. Machine Learning Research}, 8:1063--1095, 2007.

\bibitem{Tsybakov:2008:INE:1522486}
A.B. Tsybakov.
\newblock {\em Introduction to Nonparametric Estimation}.
\newblock Springer Publishing Company, Incorporated, 1st edition, 2008.

\bibitem{villmann10mathematical}
T.~Villmann and S.~Haase.
\newblock Mathematical aspects of divergence based vector quantization using
  {F}rechet-derivatives, 2010.
\newblock University of Applied Sciences Mittweida.

\bibitem{Wang-Kulkarni-Verdu2009}
Q.~Wang, S.R. Kulkarni, and S.~Verd{\'u}.
\newblock Divergence estimation for multidimensional densities via
  $k$-nearest-neighbor distances.
\newblock {\em IEEE Transactions on Information Theory}, 55(5), 2009.

\bibitem{Wolsztynski85minimum}
E.~Wolsztynski, E.~Thierry, and L.~Pronzato.
\newblock Minimum-entropy estimation in semi-parametric models.
\newblock {\em Signal Process.}, 85(5):937--949, 2005.

\bibitem{zhang11independencetest}
K.~Zhang, J.~Peters, D.~Janzing, and B.~Scholkopf.
\newblock Kernel-based conditional independence test and application in causal
  discovery.
\newblock In {\em Uncertainty in Artificial Intelligence (UAI)}, 2011.

\end{thebibliography}
}

\end{document}